\documentclass[secthm,seceqn,amsthm,ussrhead]{amsart}
\usepackage{amsmath,latexsym}
\usepackage[psamsfonts]{amssymb}
\usepackage{times}
\usepackage[mathcal]{euscript}
\numberwithin{equation}{section}

\renewcommand{\epsilon}{\varepsilon}

\newcommand{\be}{\begin{equation}}
\newcommand{\ee}{\end{equation}}

{\bf}{\it}

\begin{document}

\author {I. E. Niyozov,\, O. I. Makhmudov}

\title {The Cauchy problem for the system of the thermoelasticity
in  $E^{n}$} \maketitle

\begin{abstract}
In this paper, we considered the problem of analytical continuation
of the solution of the system equations of the thermoelasticity in
spacious bounded domain from its values and values of its strains on
part of the boundary of this domain, i.e., the Cauchy's problem.
\end{abstract}

\text{Key words}: {\small the Cauchy problem, system theory of
elasticity, elliptic system, ill-posed problem, Carleman matrix,
regularization.}

\section{Introduction}

In this paper, we considered the problem of analytical continuation
of the solution of the system equations of the thermoelasticity in
spacious bounded domain from its values and values of its strains on
part of the boundary of this domain, i.e., the Cauchy's problem.

Since, in many actual problems, either a part of the boundary is
inaccessible for measurement of displacement and tensions or only
some integral characteristic are available. In experimental study of
the stress-strain state of actual constructions, we can make
measurements only on the accessible part of the surface.

In a practical investigation of experimental dates or diagnostic
moving abject arise problems of estimation concerning deformed
position of the object. Solution of the problems by using well known
classical propositions  is connected to difficulties of absence of
experimental dates which is necessary for formulation of boundary
value (classical) conditions.

Therefore it is necessary consider the problem of continuation for
solution  of elasticity system of equations to the domain by values
of solutions and  normal derivatives in the part of boundary of
domain.

System equation of  the thermoelasticity is elliptic. Therefore the
problem Cauchy for this system is ill-posed. For ill-posed problems,
one does not prove the existence theorem: the existence is assumed a
priori. Moreover, the solution is assumed to belong to some given
subset of the function space, usually a compact one \cite{MML}. The
uniqueness of the solution follows from the general Holmgren theorem
\cite{IGP}. On establishing uniqueness in the article studio of
ill-posed problems, one comes across important questions concerning
the derivation of estimates of conditional stability and the
construction of regularizing operators.

Our aim is to construct an approximate solution using the Carleman
function method.

Let $x=(x_{1} ,...,x_{n} )$ and $y=(y_{1} ,...,y_{n} )$ be points of
the $n$-dimensional Euclidean space $E^{n} $, $D$ a bounded simply
connected  domain in $E^{n} $, with piecewise-smooth boundary
consisting of a piece   $\Sigma $ of the plane $y_{n} =0$ and a
smooth surface $S$ lying in the    half-space $y_{n} >0.$

Suppose that $n+1$-component vector function

$U(x)=(u_{1} (x),...,u_{n} (x),u_{n+1} (x))$ satisfied in $D$ the
system equations of the thermoelasticity [3]:

\begin{equation}\label{(1)}
B(\partial _{x},\omega )U(x)=0.
\end{equation}

Where
\begin{equation*}
B(\partial _{x},\omega )=\left\| B_{k\,j}(\partial _{x},\omega )
\right\| _{{(n+1)}\times{(n+1)}},
\end{equation*}
moreover
\begin{equation*}
B_{k\, j}(\partial _{x},\omega ) =\delta _{k\, j} (\mu\Delta+\rho
\omega^{2}) + (\lambda +\mu )\frac{\partial ^{2} }{\partial x_{k}
\partial x_{j} }\, ,\quad k,j=1,...,n,
\end{equation*}

\begin{equation*}
B_{k\,(n+1)}(\partial _{x},\omega ) =-\gamma \frac{\partial
}{\partial x_{(n+1)} } ,\quad k=1,...,n,
\end{equation*}

\begin{equation*}
B_{(n+1)\,j}(\partial _{x},\omega ) =i\omega\eta\frac{\partial
}{\partial x_{j} } ,\quad j=1,...,n,
\end{equation*}

\begin{equation*}
B_{(n+1)\,(n+1)}(\partial _{x},\omega )
=\Delta+\frac{i\omega}{\theta},\quad
\end{equation*}
$\delta_{ij}$ is the Kronecker delta,and where $\lambda ,\, {\kern
1pt} \mu ,\, {\kern 1pt}\rho ,\, \omega,{\kern 1pt}\theta $ is
coefficients which characterizing medium, satisfying the conditions

\begin{equation*}
\mu >0,{\kern 1pt} \, \, 3\lambda +2\mu >0,\,\, {\kern 1pt} \, \,
\rho>0,\,\,{\kern 1pt} \theta >0,\, \,{\kern
1pt}\frac{\gamma}{\eta}>0.
\end{equation*}

The system \eqref{(1)} may be write in the following way:

\begin{equation*}\label{(2)}
\left\{\begin{array}{llll} {\mu \Delta u+(\lambda +\mu)
graddivu-\gamma gradgv+\rho \omega ^{2} u =0,}  \\ { \Delta
v+\frac{i\omega}{\theta}v+ i\omega\eta divu =0,}
\end{array}\right.
\end{equation*}
where $U(x)=(u(x),v(x))$.

That system is elliptic. As,it characteristic matrix is

\begin{equation*}
 \chi(\xi)=\left\|\begin{array}{ccccc}
 (\lambda+\mu)\xi_{1}^{2}+
\mu\sum_{i=1}^{n}\xi_{i}^{2}&
(\alpha+\mu)\xi_{1}\xi_{2}&\cdots&(\alpha +\mu)\xi_{1}\xi_{n}&{0}\\
(\alpha+\mu)\xi_{2}\xi_{1}&(\lambda+\mu)\xi_{2}^{2}+\mu\sum_{i=1}^{n}\xi_{i}^{2}&\cdots
&(\alpha+ \mu)\xi_{2}\xi_{n}&{0}\\
\cdots &\cdots&\cdots&\cdots&\cdots\\
{0}&{0}&\cdots&{0}&{1}
\end{array}\right\|,
\end{equation*}
and for arbitrary real $\xi=(\xi_{1},...,\xi_{n})$ satisfying
conditions $\sum_{i=1}^{n}\xi_{i}^{2}=1,$ we have
$$
det\chi(\xi)=\mu^{2}(\lambda+\mu)\neq0
$$
\subsection*{Statement of the problem.} Find a regular solution $U$
of system \eqref{(1)} in the domain $D$ using its Cauchy data on the
surface $S$:
\begin{equation}\label{(3)}
U(y)=f(y),\quad \quad R(\partial _{y} ,\nu(y))U(y)=g(y),\quad y\in
S,
\end{equation}
where $R(\partial _{y} ,\nu(y))$ is the stress operator, i.e.,

\begin{equation*}
R(\partial_{y},\nu(y))=\left\| R_{kj}(\partial _{y}
,\nu(y),\gamma)\right\|_{(n+1)\times(n+1)}=\left\|\begin{array}{ccccc}
T & & -\gamma\nu_{1}\\
 & &-\gamma\nu_{2}\\
...&...&...&\\...&...&-\gamma\nu_{n}\\
0&0&\frac{\partial}{\partial\nu}\end{array}\right\|,
\end{equation*}
\begin{equation*}
T=T(\partial _{y} ,\nu)=\left\| T_{k\, j} (\partial _{y}
,\nu)\right\| _{n\times n},
\end{equation*}
\begin{equation*}
T_{k j} (\partial _{y},\nu)=\lambda \, \nu_{k} \frac{\partial }
{\partial y_{j} } +\mu \, \nu_{j} (y)\frac{\partial }{\partial y_{k}
} +(\mu +\lambda)\, \delta _{k\, j} \frac{\partial }{\partial
\\\nu(y)},\,\,\,\ k,j=1,...,n,
\end{equation*}
$\nu(y)=(\nu_{1} (y),...,\nu_{n} (y))$ is the unit outward normal
vector on $\partial D$ at a point $y$, $f=(f_{1} ,\ldots ,f_{n+1}
),$ $g=(g_{1} ,\ldots ,g_{n+1} )$ are given continuous vector
functions on $S.$

 \section{Construction of the matrix Carleman and approximate solution
 for the domain type's cap}

It is well known, that any regular solution $U(x)$ system
\eqref{(1)} is specified by the formula

\begin{equation}\label{(4)}
2U(x)=\int _{\partial D} \Bigg(\Psi (x-y,\omega)\{ R(\partial _{y}
,\nu(y))U(y)\} -
\end{equation}
$$
-\{{\widetilde R(\partial _{y},\nu(y))}{\widetilde\Psi
(y-x,\omega)}\} ^{*} U(y)\Bigg)ds_{y} ,{\kern 1pt} \; x\in D,
$$
where symbol $^*-$ is denote of operation transposition, $\Psi
(y,x)$ matrix of fundamental solutions system equation of
steady-state oscillations of the thermoelasticity:where
\begin{equation*} \Psi  (x,\omega)=\left\| \Psi _{k\, j}
(x,\omega)\right\| _{(n+1)\times (n+1)},
\end{equation*}

\begin{equation*}
\Psi _{k\,j}(x,\omega)=\sum _{l=1}^{3}
\Bigg[(1-\delta_{k(n+1)})(1-\delta_{j(n+1)})\left(\frac{\delta_{kj}}{2\pi\mu}\delta_{nl}
-\alpha_{l}\frac{\partial^{2}}{\partial x_k\partial x_j}\right )+
\end{equation*}
\begin{equation*}
+\beta _{l}\left(i\omega\eta\delta_{k(n+1)}(1-\delta_{j(n+1)})
 \frac{\partial}{\partial x_j}-\gamma\delta_{j(n+1)}(1-\delta_{k(n+1)})
 \frac{\partial}{\partial
 x_k}\right)+\delta_{k(n+1)}\delta_{j(n+1)}\gamma_l\Bigg]
 \frac{\exp(i\lambda_l\mid x\mid)}{\mid x\mid},
\end{equation*}

\begin{equation*}
\alpha_l=\frac{(-1)^{l}(1-i\omega\theta^{-1}\lambda^{-1}(\delta_{1l}+\delta_{2l})}
{2\pi(\lambda+2\mu)(\lambda_{2}^{2}-\lambda_{1}^{2})}-\frac{\delta_{3l}}
{2\pi\rho\omega^{2}}\,, \,\,l=1,2,3;\,\,
  \sum_{l=1}^{3}\alpha_{l}=0,
\end{equation*}
\begin{equation*}
\beta_l=\frac{(-1)^{l}(\delta_{1l}+\delta_{2l})}
{2\pi(\lambda+2\mu)(\lambda_{2}^{2}-\lambda_{1}^{2})}\,,
\,\,l=1,2,3;\,\,
  \sum_{l=1}^{3}\beta_{l}=0,
\end{equation*}
\begin{equation*}
\gamma_l=\frac{(-1)^{l}(\lambda_{l}^{2}-k_{1}^{2})(\delta_{1l}+\delta_{2l})}
{2\pi(\lambda_{2}^{2}-\lambda_{1}^{2})}\,, \,\,l=1,2,3;\,\,
  \sum_{l=1}^{3}\gamma_{l}=0,\,\,k_{1}^{2}=\rho\omega^{2}(\lambda+2\mu)^{-1},
\end{equation*}
\begin{equation*}
\widetilde\Psi(x,\omega)=\|\widetilde\Psi_{kj}(x,\omega)\|_{(n+1)\times(n+1)},
\widetilde\Psi_{kj}(x,\omega)=\Psi_{kj}(-x,\omega),
\end{equation*}
\begin{equation*}
\widetilde R(\partial _{y} ,\nu(y))=\left\|\begin{array}{ccccc}
T & & -i\omega\nu_{1}\\
 & &-i\omega\nu_{2}\\
...&...&...&\\...&...&-i\omega\nu_{n}\\
0&0&\frac{\partial}{\partial\nu}\end{array}\right\|.
\end{equation*}

\textbf{Definition.} \textit{ By the Carleman matrix of problem
\eqref{(1)},\eqref{(3)} we mean an $(n+1)\times (n+1)$ matrix $\Pi
(y,x,\omega,\tau )$ depending on the two points $y,x$ and positive
numerical number parameter $\tau $ satisfying the following two
conditions:}
$$
1){\kern 1pt} \; \Pi (y,x,\omega,\tau )=\Psi
(x-y,\omega)+G(y,x,\tau),
$$
where matrix $G(y,x,\tau )$ satisfies system \eqref{(1)} with
respect to the variable $y$ in the domain $D$, and $\Psi (y,x)$ is a
matrix of the fundamental solutions of system \eqref{(1)};

$$
2){\kern 1pt} \; \int _{\partial D\backslash S} \left(|\Pi
(y,x,\omega,\tau )| +|R(\partial _{y} ,\nu)\Pi (y,x,\omega,\tau
)|\right)ds_{y} \le \varepsilon (\tau ),
$$
where $\varepsilon (\tau )\to 0,$ as $\tau \to \infty ;$ here $|\Pi
|$ is the Euclidean norm of the matrix $\Pi =||\Pi _{i\, j}
||_{(n+1)\times(n+1)} ,$
 i.e., $|\Pi |=(\sum _{i,j=1}^{n+1} \Pi _{i\, j}^{2} )^{\frac{1}{2} } .$
 In particular, $|U|=\left(\sum _{m=1}^{n+1} u_{m}^{2})\right)^{\frac{1}{2} } .$

From the definition Karleman matrix it follows that

\textbf{Theorem 1.} \textit{ Any regular solution $U(x)$ of system
\eqref{(1)} in the domain $D$ is specified by the formula}

\begin{equation}\label{(6)}
2U(x)=\int _{\partial D} (\Pi (y,x,\omega,\tau )\{ R(\partial _{y}
,\nu)U(y)\} -
\end{equation}
$$
-\{\widetilde R(\partial _{y} ,\nu)\Pi (y,x,\omega,\tau )\} ^{*}
U(y))ds_{y}, \quad x\in D,
$$
where $\Pi (y,x,\omega,\tau )$ is matrix Carleman.

Using the matrix Carleman, easily conclude the estimate stability of
solution of the problem \eqref{(1)}, \eqref{(3)} and also indicate
effective method decision this problem as [6-8].

With a view to construct approximate solution of the problem
\eqref{(1)}, \eqref{(3)} we construct the following matrix:

\begin{equation}\label{(7)}
\Pi (y,x,\omega)=\left\| \Pi_{kj}(y,x,\omega)\right\|
_{(n+1)\times(n+1)},
\end{equation}
\begin{equation*}
\Pi_{kj}(y,x,\omega)=\sum _{l=1}^{3}
\Bigg[(1-\delta_{k(n+1)})(1-\delta_{j(n+1)})\left(\frac{\delta_{kj}}{2\pi\mu}\delta_{nl}
-\alpha_{l}\frac{\partial^{2}}{\partial x_k\partial x_j}\right )+
\end{equation*}
\begin{equation*}
+\beta _{l}\left(i\omega\eta\delta_{k(n+1)}(1-\delta_{j(n+1)})
 \frac{\partial}{\partial x_j}-\gamma\delta_{j(n+1)}(1-\delta_{k(n+1)})
 \frac{\partial}{\partial
 x_k}\right)+\delta_{k(n+1)}\delta_{j(n+1)}\gamma_l
 \Bigg]\Phi(y,x,k_{l})
\end{equation*}
where

\begin{equation}\label{(9)}
C_{n} K(x_{_{n} } )\Phi (y,x,k)=\int _{0}^{\infty }
Im[\frac{K(i\sqrt{u^{2} +s} +y_{n} )}{i\sqrt{u^{2} +s} +y_{n} -x_{n}
} ]\frac{\psi (ku)\, du}{\sqrt{u^{2} +s}},
\end{equation}

$\psi (ku)=\left\{\begin{array}{l} {uJ_{0} (ku),\quad n=2m,\quad m\ge 1,} \\
{\cos ku,\quad n=2m+1,\quad m\ge 1,} \end{array}\right. $ $J_{0}
(u)$-Bessel function of order zero,
$$
s=(y_{1} -x_{1} )^{2} +...+(y_{n-1} -x_{n-1} )^{2},\quad C_{2} =2\pi
$$

$$
C_{n} =\left\{\begin{array}{l} {(-1)^{m} \cdot 2^{-n} (n-2)\pi
\omega _{n}
 (m-2)!,\quad n=2m} \\ {(-1)^{m} \cdot 2^{-n} (n-2)\pi \omega _{n}
 (m-1)!,\quad n=2m+1.} \end{array}\right.
$$

 $K(\omega ),{\kern 1pt} \; \omega =u+iv$  ($u,{\kern 1pt} \; v$ are real),
 is an entire function taking real values on the real axis and satisfying
 the conditions $K(u)\ne \infty ,\quad \left|u\right|<\infty ,$
$K(u)\not =0,{\kern 1pt} \; \mathop{\sup }\limits_{v\ge 1} |\exp \nu
\,
 \left|Imk\right|K^{(p)} (\omega )|=\quad=M(p,u)<\infty ,{\kern 1pt} \; p=0,...,m,
 {\kern 1pt} \; u\in R^{1} .$

 In work \cite{ShYa} proved.

\textbf{Lemma 1.} \textit{ For  function $\Phi (y,x,k)$ the  formula
is valid }
\begin{equation*}\label{(10)}
{\kern 1pt} {\kern 1pt} C_{n} \Phi (y,x,k)=\varphi _{n} (ikr)+g_{n}
(y,x,k),\quad r=|y-x|,
\end{equation*}
where $\varphi _{n} $ -fudamental solution Helmholtz equation,
$g_{n} (y,x,k)$ is a regular function that is defined for all $y$
and $x$ satisfies Helmholtz equation: $\Delta (\partial _{y} )g_{n}
-k^{2} g_{n} =0.$

In \eqref{(9)} we assume the function $K(\omega )=\exp (\tau \omega
)$. Then

$$
\Phi (y,x,k)=\Phi _{\tau } (y-x,k),
$$
\begin{equation} \label{(11)}
C_{n} \Phi _{\tau } (y-x,k)=\frac{\partial ^{m-1} }{\partial s^{m-1}
} \int _{0}^{\infty } Im[\frac{\exp \, \tau (i\sqrt{u^{2} +s} +y_{n}
-x_{n} )} {i\sqrt{u^{2} +s} +y_{n} -x_{n} } ]\frac{\psi (ku)\,
du}{\sqrt{u^{2} +s} } =
\end{equation}
$$
=\exp \, \tau (y_{n} -x_{n} )\, \frac{\partial ^{m-1} }{\partial
s^{m-1} } \int _{0}^{\infty } \left[\right. -\cos \tau \sqrt{u^{2}
+\alpha ^{2} }+(y_{n} -x_{n} )\frac{\sin \tau \sqrt{u^{2} +s}
}{\sqrt{u^{2} +s} } \left.
 \right]\, \psi (ku)du,
$$
\[\Phi '_{\tau } (y-x,k)=\frac{\partial \Phi _{\tau } }{\partial \tau } .\]

\[C_{n} \Phi '_{\tau } (y-x,k)=\exp \, \tau (y_{n} -x_{n} )\,
\frac{\partial ^{m-1} }{\partial s^{m-1} } \int _{0}^{\infty }
\frac{\sin \tau \sqrt{u^{2} +s} }{\sqrt{u^{2} +s} } \psi (ku)du,\]

\begin{equation*} \label{(12)}
C_{n} \Phi '_{\tau } (y-x,k)=\exp \, \tau (y_{n} -x_{n} )\,
\frac{\partial ^{m-1} }{\partial s^{m-1} } \psi '_{\tau } (k,s),
\end{equation*}
$$
\psi '_{\tau } =\left\{\begin{array}{l} {0,\quad\qquad \tau <k} \\
{\cos \sqrt{s(\tau ^{2} -k^{2} )} ,\quad n=2m\,} \\
{\frac{1}{2} \pi \, J{}_{0}
 (\sqrt{s(\tau ^{2} -k^{2} )} ),\quad \tau >k} \end{array}\right.
 $$

Now  in formul \eqref{(7)},  and \eqref{(9)} to take $\Phi
(y,x,k)=\Phi _{\tau } (y-x,k)$, we construct matrix $\Pi
(y,x,\omega)=\Pi (y,x,\omega,\tau )$

From Lemma 1 we obtain.

\textbf{ Lemma 2.} \textit{ The matrix $\Pi (y,x,\omega,\tau )$
given by \eqref{(7)} and \eqref{(9)} is Carleman's matrix for
problem \eqref{(1)}, \eqref{(3)}.}

\textbf{ Proof.} By \eqref{(7)},  \eqref{(9)} and Lemma 1 we have
$$
\Pi (y,x,\omega,\tau )=\Psi (y,x,\omega)+G(y,x,\tau ),
$$
where
$$
G(y,x,\tau )=\left\| G_{kj}(y,x,\tau)\right\|,
$$

$$
G_{k\, j}(y,x,\tau )=\sum _{l=1}^{3}
\Bigg[(1-\delta_{k(n+1)})(1-\delta_{j(n+1)})\left(\frac{\delta_{kj}}{2\pi\mu}\delta_{nl}
-\alpha_{l}\frac{\partial^{2}}{\partial x_k\partial x_j}\right )+
$$
$$
+\beta _{l}\left(i\omega\eta\delta_{k(n+1)}(1-\delta_{j(n+1)})
 \frac{\partial}{\partial x_j}-\gamma\delta_{j(n+1)}(1-\delta_{k(n+1)})
 \frac{\partial}{\partial
 x_k}\right)+
$$$$
+\delta_{k(n+1)}\delta_{j(n+1)}\gamma_l
 \Bigg]g_{n}(y,x,k_{l} ,\tau ),\quad k,j=1,...,n+1
$$

By a straightforward calculation, we can verify that the matrix
$G(y,x,\tau )$ satisfies system \eqref{(1)} with respect to the
variable $y$ everywhere in $D$. By using \eqref{(7)}, \eqref{(9)}
and \eqref{(11)} we obtain

\begin{equation} \label{(13)}
\int _{\partial D\backslash S} \Bigg(\Big|\Pi (y,x,\omega,\tau
)\Big|+ \Big|R(\partial _{y} ,\nu)\Pi (y,x,\omega,\tau )\Big|\Bigg),
ds_{y} \le C_{1} (x) \tau ^{m} \exp (-\tau \, x_{n} ),
\end{equation}

where $C_{1}(x)$ some bounded function inside $D.$ The lemma is
thereby proved.

Let us set
\begin{equation} \label{(14)}
2U_{\tau } (x)=\int _{S} [\Pi (y,x,\omega,\tau )\{ R(\partial _{y}
,\nu)U(y)\} - \{\widetilde R(\partial _{y} ,n)\Pi (y,x,\omega,\tau
)\} ^{*} U(y)]ds_{y}.
\end{equation}

The following theorem holds.

\textbf{Theorem 1.} \textit{ Let $U(x)$ be a regular solution of
system \eqref{(1)} in $D$ such that }

\begin{equation} \label{(15)}
|U(y)|+|R(\partial _{y} ,\nu)U(y)|\le M,{\kern 1pt} \; y\in
\partial D\backslash S.
\end{equation}
Then for $\tau \ge 1$ the following estimate is valid:

\[|U(y)-U_{\tau } (y)|\le MC_{2} (x)\tau ^{m} \exp (-\tau \, x_{n} ).\]

\textbf{Proof.} By formula \eqref{(6)} and \eqref{(14)}, we have

$$
2|U(x)-U_{\tau } (x)|=\int _{\partial D\backslash S} [\Pi
(y,x,\omega,\tau ) \{ R(\partial_{y},\nu)U(y)\} -\{\widetilde
R(\partial _{y},\nu)\Pi (y,x,\omega,\tau )\} ^{*}
 U(y)]ds_{y} .
$$
Now on the basis of \eqref{(13)} and \eqref{(15)} we obtain required
inequality. The theorem is thereby proved.

Now we write out a result that allows us to calculate $U(x)$
approximately if, instead of $U(y)$ and $R(\partial _{y} ,\nu)U(y),$
their continuous approximations $f_{\delta } (y)$ and $g_{\delta }
(y)$ are given on the surface $S$:

\begin{equation} \label{(16)}
\mathop{\max }\limits_{S} |f(y)-f_{\delta } (y)|+\mathop{\max }
\limits_{S} |R(\partial _{y},\nu)U(y)-g_{\delta } (y)|\le \delta ,
{\kern 1pt} \; 0<\delta <1.
\end{equation}

We define a function $U_{\tau \, \delta } (x)$ by setting

\begin{equation} \label{(17)}
2U_{\tau \, \delta } (x)=\int _{S} [\Pi (y,x,\omega,\tau )g_{\delta
} (y)-\{ \widetilde R(\partial _{y} ,\nu)\Pi (y,x,\omega,\tau )\}
^{*} f_{\delta }(y)]ds_{y},
\end{equation}
where

$$
\tau =\frac{1}{x_{n}^{0} } ln\frac{M}{\delta } ,\quad x_{n}^{0}
=\mathop{\max }\limits_{D} x_{n} ,\quad x_{n} >0.
$$

\textbf{Theorem 2.} \textit{ Let $U(x)$ be a regular solution of
system \eqref{(1)} in $D$ satisfying condition \eqref{(15)}. Then
the following estimate is valid:}

$$
|U(x)-U_{\tau \, \delta } (x)|\le C_{3} (x)\delta ^{\frac{x_{n} }
{x_{n}^{0} } } \left(ln\frac{M}{\delta } \right)^{m} ,\quad x\in D.
$$

\textbf{Proof.} From formula \eqref{(6)} and \eqref{(17)}
 we have
$$
2(U(x)-U_{\tau \, \delta } (x))=\int _{S} [\Pi (y,x,\omega,\tau )\{
R(\partial _{y} ,\nu) U(y)-g_{\delta } (y)\} -
$$
$$
-\{\widetilde R(\partial _{y} ,\nu)\Pi (y,x,\omega,\tau )\}^{*}
(U(y)-f_{\delta } (y))]ds_{y} +
$$
$$
+\int _{\partial D\backslash S} [\Pi (y,x,\omega,\tau )\{ R(\partial
_{y} ,\nu) U(y)\} -\{\widetilde R(\partial _{y} ,\nu)\Pi
(y,x,\omega,\tau )\} ^{*} U(y)]ds_{y}.
$$
By the assumption of the theorem and inequalities \eqref{(13)},
\eqref{(15)} and \eqref{(16)} for the any $x\in D,$ we obtain

$$
|U(x)-U_{\tau \, \delta } (x)|=C_{2} (x)\delta \tau ^{m} \exp
 \tau (x_{n}^{0} -x_{n} )+C_{1} (x)\tau ^{m} \exp (-\tau \, x_{n} )\le
$$
$$
\le C_{3} (x)\tau ^{m} (M+\delta \exp \tau \, x_{n}^{0} )\exp (-\tau
\, x_{n} ).
$$
Now, it to take $\tau =\frac{1}{x_{n}^{0} } ln\frac{M}{\delta } ,$
 then we obtain to proof theorem. The theorem is thereby proved.

\textbf{Theorem 3.} \textit{ Let $U(x)$ be a regular solution of
system \eqref{(1)} in $D$ satisfying conditions }

$$
|U(y)|+|R(\partial _{y} ,\nu)U(y)|\le M,{\kern 1pt} \; y\in \partial
D\backslash S,
$$
$$
|U(y)|+|R(\partial _{y} ,\nu)U(y)|\le \delta ,\quad 0<\delta
<1,\quad y\in S.
$$
Then
$$
|U(x)|\le C_{4} (x)\delta ^{\frac{x_{n} }{x_{n}^{0} } }
\left(ln(\frac{M}{\delta } )\right)^{m},
$$
where $C_{4} (x)=\tilde{C}\int _{\partial D} \frac{1}{r^{n} }
 ds_{y} ,\quad \tilde{C}-$ constant depending on $\lambda ,
 \mu ,\omega .$

\textbf{Proof.} On the basis of Theorem 2 we obtain

$$
|U(x)|\le |U_{\tau } (x)|+MC_{2} (x)\tau ^{m} \exp (-\tau \, x_{n}
).
$$
Next from the condition theorem and \eqref{(7)}, \eqref{(9)} we
obtain
$$
|2U_{\tau \, \delta } (x)|=\left|\int _{S} [\Pi (y,x,\omega,\tau )
\{ R(\partial _{y} ,\nu)U(y)\} -\{ \widetilde R(\partial _{y}
,\nu)\Pi (y,x,\omega,\tau )\}^{*} (U(y)\left. )\right]ds_{y}
\right|\le
$$
$$
\le \int _{S} \left(|\Pi (y,x,\omega,\tau )|+|R(\partial _{y} ,\nu)
\Pi (y,x,\omega,\tau )|\right)\,\left(|U(y)|+|R(\partial _{y}
,\nu)U(y)|\right)ds_{y} \le
$$
$$
\le\delta \int_{S}\left(|\Pi (y,x,\omega,\tau)|+|R(\partial_{y},\nu)
\Pi(y,x,\omega,\tau )|\right)ds_{y} \le C_{3} (x)\delta \tau ^{m}
\exp (\tau \, x_{n}^{0} -\tau \, x_{n} ).
$$
Then
$$
|U(x)|\le C_{4} (x)\tau ^{m} \exp (-\tau \, x_{n} )(M+\delta \exp
\tau \, x_{n}^{0} ).
$$
Next if we take $\tau =\frac{1}{x_{n}^{0} } ln\frac{M}{\delta } ,$
then we obtain stability estimate:

$$|U(x)|\le C_{4} (x)\delta ^{\frac{x_{n} }{x_{n}^{0} } }
\left(ln(\frac{M}{\delta } )\right)^{m}.
$$

The theorem is thereby proved.

From proved above theorems we obtain

\textbf{Corollary 1.} \textit{ The limit relation }

$$
\mathop{\lim }\limits_{\tau \to \infty } U_{\tau } (x)=U(x), {\kern
1pt} \quad \mathop{\lim } \limits_{\delta \to 0} U_{\tau \, \delta }
(x)=U(x)
$$
hold uniformly on each compact subset of $D.$

 \section{ Regularization of solution of the problem \eqref{(1)},
  \eqref{(3)} for a domain of cone type}

Let $x=(x_{1} ,\ldots,x_{n} )$ and $y=(y_{1} ,\ldots ,y_{n} )$ be
points in $E^{n} ,\, \, {\kern 1pt} {\kern 1pt} D_{\rho } $ be a
bounded simply connected domain in $E^{n} $ whose boundary consists
of a cone surface

\[\Sigma :\quad \alpha _{1} =\tau _{\rho } \, y_{n} ,{\kern 1pt} \;
\alpha _{1}^{2} =y_{1}^{2} +\ldots +y_{n-1}^{2} ,{\kern 1pt} \; \tau
_{\rho } =tg\frac{\pi }{2\rho } ,{\kern 1pt} \; y_{n} >0,{\kern 1pt}
\; \rho >1\] and a smooth surface $S,$ lying in the cone. Assume
$x_{0} = (0,...0,x_{n} )\in D_{\rho } $.

We constract Karleman matrix. In formula \eqref{(7)}, \eqref{(9)} to
take

\[K(\omega )=E_{\rho } \left[\tau \left(\omega -x_{n} \right)\right],
\, \, \, \, \, \, \, \tau >0,\, \, \, \, \rho >1.\, \] Then

\[\Phi (y,x,k)=\Phi _{\tau } (y-x,k),  k>0\]

\begin{equation} \label{(18)}
C_{n} \Phi _{\tau } (y-x,k)=\frac{\partial ^{m-1} }{\partial s^{m-1}
} \int _{0}^{\infty } Im\, [\frac{E_{\rho } \left(\tau (i\sqrt{u^{2}
+s} +y_{n} -x_{n} )\right)}{i\sqrt{u^{2} +s} +y_{n} -x_{n} } ]
\frac{\psi (ku)\, du}{\sqrt{u^{2} +s} }
\end{equation}

\[\Phi '_{\tau } (y-x,k)=\frac{\partial \Phi _{\tau } }{\partial \tau } .\]

\begin{equation*} \label{(19)}
C_{n} \Phi '_{\tau } (y-x,k)=\, \frac{\partial ^{m-1} }{\partial
s^{m-1} }
 \int _{0}^{\infty } Im\left\{\, E'_{\rho } \, \left[\tau
 (i\sqrt{u^{2} +s} +y_{n} -x_{n} \right]\, \, \right\}\frac{\psi (ku)\,
 du}{\sqrt{u^{2} +s} },
 \end{equation*}
where $\quad E_{\rho } (w)-$ Mittag-L\"{o}ffer`s a entire function
\cite{MMD}. For the functions $\Phi _{\tau } (y-x,k)$ holds Lemma 1
and Lemma 2.

Now again to denote by $U_{\tau } (x), {\kern 1pt} \, \, \, \,
U_{\tau \, \delta } (x)$ as \eqref{(14)} and \eqref{(17)}. Then
holds analogical theorem as \textbf{Theorem} \textbf{1,2,3.}

For    $n=3$ we reduce entirely.

Suppose that $D_{\rho } $ is bounded simple connected domain from
$E^{3} $ with boundary consisting of part $\sum  $ of the surface of
the cone

\[y_{1}^{2} +y_{2}^{2} =\tau _{\rho } \, y_{3}^{2} ,\quad \tau _{\rho } =
tg\frac{\pi }{2\rho } ,\quad \rho >1,\quad y_{3} >0,\] and of a
smooth portion of the surface $S$ lying inside the cone. Assume
$x_{0} =(0,0,x_{3} )\in D_{\rho } $.

We construct Carleman`s matrix. In formula \eqref{(7)}, \eqref{(9)}
we take

\begin{equation} \label{(20)}
\Phi _{\tau } (y,x,k)=\frac{1}{4\pi ^{2} E_{\rho } (\tau
^{\frac{1}{\rho } } x_{3} )} \int _{0}^{\infty } Im\frac{E_{\rho }
(\tau ^{\frac{1}{\rho } } w)} {i\sqrt{u^{2} +s} +y_{3} -x_{3} }
\frac{\cos ku\, du}{\sqrt{u^{2} +s} },
\end{equation}
where $w=i\sqrt{u^{2} +s} +y_{3}$. For the functions $\Phi _{\tau }
(y,x,k)$ holds Lemma 1.

If follows from the properties of $E_{\rho } (w)$ that for $y\in
\Sigma,\quad 0<u<\infty $ the function $\Phi _{\tau } (y,x,k)$
defined by \eqref{(18)} its gradient and second partial derivatives

\[\frac{\partial ^{2}\Phi _{\tau }(y,x,k)}{\partial y_{k} \partial y_{j} }
,\quad k,j=1,2,3,\] tend to zero as $\tau \to \infty $ for a fixed
$x\in D_{\rho } .$

Then from \eqref{(7)} we find that the matrix $\Pi (y,x,\omega,\tau
)$ and its stresses $R(\partial _{y} ,\nu)\Pi (y,x,\omega,\tau )$
also tend to zero as $\tau \to \infty $ on $y\in \Sigma ,$ i.e.,
$\Pi (y,x,\omega,\tau )-$ is the Carleman matrix for the domain
$D_{\rho } $ and the part $\Sigma $ of the boundary.

For the $U(x)-$ regular solution system \eqref{(1)} following
integral formula holds

\begin{equation*} \label{(21)}
2U(x)=\int _{\partial D_{\rho } } [\Pi (y,x,\omega,\tau )\{
R(\partial _{y} ,\nu)U(y)\} -\{\widetilde R(\partial _{y} ,\nu)\Pi
(y,x,\omega,\tau )\} ^{*} U(y)]ds_{y}.
\end{equation*}

By $x\in D_{\rho } $ we denote $U_{\tau } (x)$ follows:

\begin{equation} \label{(22)}
2U_{\tau } (x)=\int _{S} [\Pi (y,x,\omega,\tau )\{ R(\partial _{y}
,\nu)U(y)\}- \{ \widetilde R(\partial _{y} ,\nu)\Pi (y,x,\omega,\tau
)\} ^{*} U(y)]ds_{y}.
\end{equation}

The following theorem holds.

\textbf{Theorem 4.} \textit{ Let $U(x)$ be a regular solution of
system \eqref{(1)} in $D_{\rho } $ such that }
\begin{equation} \label{(23)}
|U(y)|+|R(\partial _{y} ,\nu)U(y)|\le M,{\kern 1pt} \; y\in \Sigma.
\end{equation}

Then for $\tau \ge 1$ the following estimate is valid:

\[|U(x_{0} )-U_{\tau } (x_{0} )|\le MC_{\rho } (x_{0} )\tau ^{3}
\exp(-\tau \, x_{3}^{\rho } ),\]

where $x_{0} =(0,0,x_{3} )\in D_{\rho } ,\quad x_{3} >0,\quad $

\[C_{\rho }(x_{0} )=C_{\rho }\int _{\Sigma }\frac{1}{r_{0}^{3}}ds_{y},
\quad r_{0} =|y-x_{0} |,\quad C_{\rho } -constant.\]

\textbf{Proof.} By analogy with proved Theorem 2 and Theorem 3 from
\eqref{(22)} and \eqref{(23)} we obtain

\[|U(x_{0} )-U_{\tau } (x_{0} {\rm \; \; }|\le M\int _{\Sigma }
[|\Pi (y,x_{0} ,\omega,\tau )|+|\widetilde R(\partial _{y} ,\nu)\Pi
(y,x_{0} ,\omega,\tau )|]ds_{y}.\]

By formula \eqref{(20)} we have following inequality:

\[\left|\Phi _{\tau } (y,x,k)\right|\le C_{\rho }^{(1)} E_{\rho }^{-1}
(\tau ^{\frac{1}{\rho } } x_{3} )r^{-1} ,\]

\[\left|\frac{\partial \Phi _{\tau } (y,x,k)}{\partial y_{i} }
\right|\le C_{\rho }^{(2)} \, \tau \, E_{\rho }^{-1} (\tau
^{\frac{1}{\rho } }x_{3} )r^{-2} \]

\[\left|\frac{\partial ^{2} \Phi _{\tau } (y,x,k)}{\partial y_{k} \partial y_{j} }
\right|\le C_{\rho }^{(3)} \tau ^{2} E_{\rho }^{-1}
 (\tau ^{\frac{1}{\rho } } x_{3} )r^{-3} .\]

Then from \eqref{(7)}

\[|\Pi (y,x,\omega,\tau )|\le C_{\rho }^{(4)} \tau ^{2} E_{\rho }^{-1}
(\tau ^{\frac{1}{\rho } } x_{3} )r^{-3} ,\]

\[|\widetilde R(\partial _{y} ,\nu)\Pi (y,x,\omega,\tau )|\le C_{\rho }^{(5)}
\tau ^{3} E_{\rho }^{-1} (\tau ^{\frac{1}{\rho } } x_{3} )r^{-4} .\]

Therefore  we obtain

\[|U(x_{0} )-U_{\tau } (x_{0} )|\le MC_{\rho } (x_{0} )\tau ^{3}
 \exp (-\tau \, x_{3}^{\rho } ),\]

where

\[C_{\rho }(x_{0} )=C_{\rho }\int _{\Sigma }\frac{1}{r_{0}^{3} }ds_{y},
\quad r_{0} =|y-x_{0}|,\quad C_{\rho } -constant.\]

The theorem is thereby proved.

Suppose that instead of $U(y)$ and $R(\partial _{y} ,\nu)U(y)$ gives
their continuous approximations $f_{\delta } (y)$ and $g_{\delta }
(y)$ such that

\[\mathop{\max }\limits_{S} |U(y)-f_{\delta } (y)|+\mathop{\max }\limits_{S}
|T(\partial _{y} ,n)U(y)-g_{\delta } (y)|\le \delta,{\kern 1pt} \;
\, \, \, 0<\delta <1.\]

Define the function $U_{\tau \, \delta } (x)$ by

\[2U_{\tau \, \delta } (x)=\int _{S} [\Pi (y,x,\omega,\tau )g_{\delta }
 (y)-\{ \widetilde R(\partial _{y} ,\nu)\Pi (y,x,\omega,\tau )\} ^{*} f_{\delta }
 (y)]ds_{y} ,\, \, \, \, \, \, \, \, \, \, \, \, \, \]

The following theorem holds

\textbf{Theorem 5.} \textit{ Let $U(x)$ is a regular solution of
system \eqref{(1)} in the domain $D_{\rho } $ satisfying the
condition \eqref{(23)}, then }

\[|U(x_{0} )-U_{\tau \, \delta } (x_{0} )|\le C_{\rho } (x_{0} )
\delta ^{q} (ln\frac{M}{\delta } )^{3} ,\]

where $\tau =(\tau _{\rho } \, R)^{-\rho } ln\frac{M}{\delta },\quad
R^{\rho } =\mathop{\max }\limits_{S} Re(i\sqrt{s} +y_{3} )^{\rho }
,$

\[q=(\frac{x_{3} }{R} )^{\rho },\quad C_{\rho } (x_{0} )=C_{\rho }
\int_{\Sigma } \left[\frac{1}{r_{0}^{3} } +\frac{1}{r_{0}^{4} }
\right]ds_{y} .\]

 The proof theorem is similar to those of Theorem 3 and 4.

\textbf{Corollary 2.} \textit{ The limit relation }

\[\mathop{\lim }\limits_{\tau \to \infty } U_{\tau } (x)=U(x),{\kern 1pt}
\, \, \, \, \, \, \, \, \, \, \, \, \, \, \, \, \mathop{\lim }
\limits_{\delta \to 0} U_{\tau \, \delta } (x)=U(x)\]

hold uniformly on each compact subset of $D_{\rho } .$


\begin{thebibliography}{11}



\bibitem{MML} M.M.Lavrent'ev. Some Ill-Posed Problems of Mathematical
Physics [in Russian], Computer Center of the Siberian Division of
the Russian Academy of Sciences, Novosibirck (1962) 92p.

\bibitem{IGP} I.G.Petrovskii. Lectures on Partial Differential Equations
[in Russians],Fizmatgiz,Moscow,(1961).

\bibitem{VDK} V.D.Kupradze, T.V.Burchuladze, T.G.Gegeliya, ot.ab.
Three-Dimensional Problems of the Mathematical Theory of Elasticity
and ... [in Russian],Nauka,Moscow,1976.

\bibitem{ShYa} Sh.Ya.Yarmukhamedov. Dokl.Acad.Nauk SSSR [Soviet Math.Dokl.],
 V.357, No.3, p.320-323.(1997).

\bibitem{MMD} M.M. Dzharbashyan, Integral Transformations and Representations
 of Functions in a Complex Domain [in Russian],Nauka, Moscow.1966.
\bibitem{MNT} O. I. Makhmudov, I. E. Niyozov, N.Tarkhanov. The Cauchy Problem
of Couple-Stress Elasticity. Contemporary Mathematics.AMS, V455,
2008.pp 207-310.

\bibitem{MN}  O. I. Makhmudov, I. E. Niyozov. The Cauchy problem for the
Lame system in infinite domains in $R^{m} $. Journal of inverse and
Ill-Posed Problems.V14.N9.2006.pp.905-924(20).
\bibitem{MN} O. I. Makhmudov, I. E. Niyozov.Regularization of a solution to the
Cauchy Problem for the System of Thermoelasticity.Contemporary
Mathematics.AMS, Primary V382, 2005,74F05, 35Q72.
\end{thebibliography}
\end{document}